Some Results on the Counterfeit Coins Problem  II

Li An-Ping

Beijing 100085, P.R.China
apli0001@sina.com

Abstract

In this paper, we will continue the earlier discussion to show some estimations to $g_1(n|_m)$ for general $n$ and $m$.

Keywords: combinatorial search, counterfeit coins problem, information theoretic bound.

In this paper, we will proceed the discussion in [1] to estimate $g_1(n|_k)$ for general $n$ and $k$. At first, we introduce some notations.

Let $\mathcal{S}$ be same as defined in the earlier paper [1], for any point $X \in \mathcal{S}$, $X = (x_1, \cdots, x_m)$, denoted by $\tilde{X} = \{x_1, \cdots, x_m\}$. Suppose that $\mathcal{D}$ is a subset, $\mathcal{D} \subseteq \mathcal{S}$, if for any point $Z \in \mathcal{D}$, there is an element $x_Z \in \tilde{Z}$ such that $x_Z \notin \tilde{X}, \forall X \in \mathcal{D}, X \neq Z$, then we call that the set $\mathcal{D}$ is $1-representable$ set, and the set $\{x_Z | Z \in \mathcal{D}\}$ is a representative set of $\mathcal{D}$. Similarly, if $\forall Z \in \mathcal{D}$, there are two elements $x_{z_1}, x_{z_2} \in \tilde{Z}$, such that $x_{z_1}$, or $x_{z_2} \notin \tilde{X}, \forall X \in \mathcal{D}, X \neq Z$, then set $\mathcal{D}$ will be called $2-representable$, and the set $\{(x_{z_1}, x_{z_2}) | Z \in \mathcal{D}\}$ is called a 2-representative set of set $\mathcal{D}$.

$S_1, \cdots, S_m$ are the $m$ sets of coins same as in [1], $|S_i| = n_i, Ct(S_i) = 1, 1 \leq i \leq m$, and

$$\mathcal{S} = \prod_{i=1}^{m} S_i,$$ $\mathcal{S}_\alpha$ represents the objective set on the direction $\alpha$ in a given algorithm. Denoted by $\Omega_k = \{\mathcal{S}_\alpha \| \alpha | = k, \alpha \in \mathcal{A}\}$, and in this paper we adopt the expression $\mathcal{S} \xrightarrow{k} \Omega_k$, or more clearly,

$$(n_1, \cdots, n_m) \xrightarrow{k} \Omega_k. \qquad (1)$$

The cases we are interested are that $\Omega_k$ contains only one or a few types of objective sets with simpler construction, since this will help to simplify searching algorithm. For instance, the objective sets in $\Omega_k$ all are $1-representable$ sets with no more than $r$ coins, in this case we will represent it as

$$(n_1, \cdots, n_m) \xrightarrow{k} (r). \qquad (2)$$

Similarly, each objective set in $\Omega_k$ is 2-representable, and the 2-representive set has the form

$$\{x_i \cdot \Delta_i | x_i \in S_j, \Delta_i \subseteq \bigcup_v S_v, |\Delta_i| \leq t, 1 \leq i \leq s\},$$

then we will write this fact as

$$(n_1, \cdots, n_m) \xrightarrow{k} (r \times s). \qquad (3)$$

More general, if $\Omega_k$ includes more than one type of objective sets, for example, two types $type1$ and $type2$, then we write it as

$$(n_1,\cdots,n_m)\xrightarrow{k} type1 \vee type2 ,\tag{4}$$

If there is an algorithm such that

$$(n_1,\cdots,n_m)\xrightarrow{i}(r) \text{ and } (n_1,\cdots,n_m)\xrightarrow{j}(s),$$

Then we will represent this case as

$$(n_1,\cdots,n_m)\xrightarrow{i,j}(r,s).\tag{5}$$

In the earlier discussion [1], [2], we have known that

$$g_1(2\mid_3) = 2,\tag{6}$$
$$g_1(5\mid_2) = 3,\tag{7}$$
$$g_1(7\mid_5) = 9,\tag{8}$$
$$g_1(11\mid_5) = 11,\tag{9}$$
$$g_1(13\mid_5) = 12,\tag{10}$$
$$g_1(17\mid_5) = 13,\tag{11}$$
$$g_1(20\mid_4) = 11,\tag{12}$$
$$g_1(5,16) = 4,\tag{13}$$

In the next, we will show more of such individual results in order to present a estimate for $g_1(n\mid_k)$ for general $n$ and $k$.

**Lemma 1**

$$g_1(4,4,5) = 4,\tag{14}$$
$$g_1(8,10) = 4,\tag{15}$$
$$g_1(2,4,10) = 4,\tag{16}$$
$$g_1(11,22) = 5,\tag{17}$$
$$g_1(7,7,7,2) = 6,\tag{18}$$
$$g_1(7,7,14) = 6,\tag{19}$$
$$g_1(8,8,8,4) = 7,\tag{20}$$
$$g_1(8,8,8,2,2) = 7,\tag{21}$$
$$g_1(35\mid_4) = 13,\tag{22}$$
$$g_1(32\mid_6) = 19,\tag{23}$$
$$g_1(46\mid_2) = 7,\tag{24}$$
$$g_1(2,19,19) = 6,\tag{25}$$
$$g_1(4,4,4,4,76) = 9,\tag{26}$$

**Lemma 2**

$$(4,4,4,4,4)\xrightarrow{5}(5),\tag{27}$$

$$(10,10) \xrightarrow{3} (4), \qquad (28)$$

$$(13,13,13,13) \xrightarrow{8} (5), \qquad (29)$$

$$(16,16) \xrightarrow{3} (10), \qquad (30)$$

$$(19,19) \xrightarrow{3} (14) \vee (2\times 7), \qquad (31)$$

$$(23,23) \xrightarrow{3} (4\times 5), \qquad (32)$$

$$(25,25) \xrightarrow{4} (8), \qquad (33)$$

$$(28,28) \xrightarrow{4} (10), \qquad (34)$$

$$(29,29) \xrightarrow{4} (11), \qquad (35)$$

$$(38,38) \xrightarrow{6} (2), \qquad (36)$$

$$(14,42) \xrightarrow{3} (22), \qquad (37)$$

$$(49,49) \xrightarrow{5} (10), \qquad (38)$$

$$(52,52) \xrightarrow{4} (35), \qquad (39)$$

$$(56,56) \xrightarrow{4,6} (40,4), \qquad (40)$$

$$(58,58) \xrightarrow{5} (14), \qquad (41)$$

$$(61,61) \xrightarrow{5} (2\times 23) \qquad (42)$$

$$(73,73) \xrightarrow{5} (22) \qquad (43)$$

The algorithms for these formulas are described in the appendix in the end of the paper.

From the formulas above, it is easy to obtain

**Corollary 1**

$$g_1(4\,|_7) = 9, \qquad (44)$$
$$g_1(8\,|_{11}) = 21, \qquad (45)$$
$$g_1(10\,|_8) = 17, \qquad (46)$$
$$g_1(13\,|_8) = 19, \qquad (47)$$
$$g_1(14\,|_9) = 22, \qquad (48)$$
$$g_1(16\,|_9) = 23, \qquad (49)$$
$$g_1(19\,|_7) = 19, \qquad (50)$$
$$g_1(22\,|_6) = 17, \qquad (51)$$
$$g_1(23\,|_8) = 23, \qquad (52)$$
$$g_1(28\,|_{14}) = 43, \qquad (53)$$
$$g_1(38\,|_3) = 10, \qquad (54)$$
$$g_1(40\,|_8) = 27, \qquad (55)$$
$$g_1(42\,|_{12}) = 41, \qquad (56)$$
$$g_1(44\,|_{15}) = 52, \qquad (57)$$
$$g_1(49\,|_5) = 18, \qquad (58)$$

$$g_1(52|_8) = 29, \tag{59}$$
$$g_1(56|_{10}) = 37, \tag{60}$$
$$g_1(58|_7) = 26, \tag{61}$$
$$g_1(61|_4) = 15, \tag{62}$$
$$g_1(70|_{10}) = 39. \tag{63}$$
$$g_1(73|_{12}) = 47, \tag{64}$$
$$g_1(76|_{19}) = 75, \tag{65}$$

Proof. The equation (43) is followed from (14) and (25). The equation (44) is easy followed from (20) and (21). (45) is followed from (6), (16) and (26). (46) is followed from (27) and (7). (47) is followed (6) and (19). The rest equations above may be obtained similarly from the equations in Lemma 1, 2, which are omitted. □

Up to now, for each positive integer $n \leq 81$, we have found a corresponding positive integer $k_0(n)$ and an algorithm for getting $g_1(n|_{k_0(n)})$. Denoted by $\mu(n) = g_1(n|_{k_0(n)}) / k_0(n)$, there is

**Proposition 1.** For each positive integer $n \leq 81$, and arbitrary positive integer $k$, then

$$\lceil k \cdot \log_3 n \rceil \leq g_1(n|_k) \leq \lceil k \cdot \mu(n) \rceil. \tag{65}$$

Proof. The left-hand side of (65) is the information theoretical bound of $g_1(n|_k)$, so we are provided to show the inequality of the right-hand hand. Clearly, it is true for the cases $k_0(n) | k$, hence it suffices to verify it for $k < k_0(n)$, these algorithms usually are not more difficult than the ones for $g_1(n|_{k_0(n)})$ and are settled easily by applying the known algorithms above, so are omitted. □

In the next is a list for $\mu(n)$ and $\log_3(n)$ for $n \leq 81$ in order to compare them each other.

| $n$ | $\mu(n)$ | $\log_3 n$ | $n$ | $\mu(n)$ | $\log_3 n$ | $n$ | $\mu(n)$ | $\log_3 n$ |
|---|---|---|---|---|---|---|---|---|
| 1 | 0 | 0 | 31 | 19/6 | 3.126 | 61 | 15/4 | 3.742 |
| 2 | 2/3 | 0.631 | 32 | 19/6 | 3.155 | 62 | 19/5 | 3.757 |
| 3 | 1 | 1 | 33 | 16/5 | 3.183 | 63 | 19/5 | 3.771 |
| 4 | 9/7 | 1.262 | 34 | 13/4 | 3.210 | 64 | 23/6 | 3.786 |

| 5 | 3/2 | 1.465 | 35 | 13/4 | 3.236 | 65 | 23/6 | 3.800 |
| --- | --- | --- | --- | --- | --- | --- | --- | --- |
| 6 | 5/3 | 1.631 | 36 | 23/7 | 3.262 | 66 | 23/6 | 3.814 |
| 7 | 9/5 | 1.771 | 37 | 10/3 | 3.287 | 67 | 31/8 | 3.827 |
| 8 | 21/11 | 1.891 | 38 | 10/3 | 3.311 | 68 | 31/8 | 3.841 |
| 9 | 2 | 2 | 39 | 27/8 | 3.335 | 69 | 31/8 | 3.854 |
| 10 | 17/8 | 2.096 | 40 | 27/8 | 3.358 | 70 | 39/10 | 3.867 |
| 11 | 11/5 | 2.183 | 41 | 31/9 | 3.380 | 71 | 43/11 | 3.880 |
| 12 | 16/7 | 2.262 | 42 | 31/9 | 3.402 | 72 | 43/11 | 3.893 |
| 13 | 19/8 | 2.335 | 43 | 7/2 | 3.424 | 73 | 47/12 | 3.905 |
| 14 | 22/9 | 2.402 | 44 | 7/2 | 3.445 | 74 | 75/19 | 3.918 |
| 15 | 5/2 | 2.465 | 45 | 7/2 | 3.465 | 75 | 75/19 | 3.930 |
| 16 | 23/9 | 2.524 | 46 | 7/2 | 3.485 | 76 | 75/19 | 3.942 |
| 17 | 13/5 | 2.579 | 47 | 32/9 | 3.505 | 77 | 4 | 3.954 |
| 18 | 8/3 | 2.631 | 48 | 32/9 | 3.524 | 78 | 4 | …… |
| 19 | 19/7 | 2.680 | 49 | 32/9 | 3.543 | 79 | 4 | …… |
| 20 | 11/4 | 2.727 | 50 | 18/5 | 3.561 | 80 | 4 | …… |
| 21 | 14/5 | 2.771 | 51 | 18/5 | 3.579 | 81 | 4 | 4 |
| 22 | 17/6 | 2.814 | 52 | 29/8 | 3.597 | | | |
| 23 | 23/8 | 2.854 | 53 | 11/3 | 3.614 | | | |
| 24 | 32/11 | 2.893 | 54 | 11/3 | 3.631 | | | |
| 25 | 3 | 2.93 | 55 | 37/10 | 3.648 | | | |
| 26 | 3 | 2.96 | 56 | 37/10 | 3.664 | | | |
| 27 | 3 | 3 | 57 | 26/7 | 3.680 | | | |
| 28 | 28/9 | 3.033 | 58 | 26/7 | 3.696 | | | |
| 29 | 31/10 | 3.065 | 59 | 15/4 | 3.712 | | | |
| 30 | 25/8 | 3.096 | 60 | 15/4 | 3.727 | | | |

List 1

From the list above we have seen that the values $\mu(n)$ are much closed to $\log_3 n$ for $n \leq 81$. The following is a result for general $n$ and $k$.

**Proposition 2**

$$\lceil k \cdot \log_3 n \rceil \leq g_1(n\vert_k) \leq \lceil k \cdot (\log_3 n + 0.076) \rceil \tag{66}$$

Proof. The left-hand side of inequality (66) is the information theoretic bound of $g_1(n\vert_k)$, we are provided to show the right-hand of the inequality. By Proposition 1, we know that it is true for $n \leq 81$, so assume $n > 81$.

Let $S_i, 1 \leq i \leq k$, same as in the paper [1], be $k$ sets of coins, $\vert S_i \vert = n, Ct(S_i) = 1, 1 \leq i \leq k$.

Suppose that $n = \lambda \cdot 3^l, 1 \leq \lambda \leq 3$, $l$ is a non-negative integer, by the assumption, $l > 3$. Let $\{\lambda_i\}_0^{25}$, $0 \leq i \leq 25$, $1 = \lambda_0 < \lambda_1 < \cdots < \lambda_{25} = 3$, which are

$1, 28/27, 10/9, 32/27, 11/9, 35/27, 4/3, 38/27, 40/27, 14/9, 44/27, 46/27, 49/27,$
$17/9, 52/27, 2, 56/27, 58/27, 61/27, 7/3, 22/9, 23/9, 70/27, 8/3, 76/27, 3.$

Suppose that $\lambda \in (\lambda_{j-1}, \lambda_j]$ for some $j$, $1 \leq j \leq 25$. For each $S_i$, $1 \leq i \leq k$, we spend $l-3$ weighings, then we will find a subset $A_i \subset S_i$, with that $|A_i| \leq \lceil \lambda \cdot 3^3 \rceil$, $Ct(A_i) = 1, 1 \leq i \leq k$.

Without loss generality, we may assume that $|A_i| = d = \lceil \lambda \cdot 3^3 \rceil, 1 \leq i \leq k$. Denoted by $d_i = \lambda_i \cdot 3^3$, $1 \leq i \leq 25$. Obviously, $d_{j-1} < d \leq d_j$, so by Proposition 1,

$$g_1(d|_k) \leq g_1(d_j|_k) \leq \lceil k \cdot \mu(d_j) \rceil.$$

Moreover, from List 1 we know

$$\max_i \{\mu(d_i) - \log_3(d_{i-1} + 1)\} < 0.076$$

Therefore,

$$g_1(d|_k) \leq \lceil k \cdot (\log_3(d_{j-1} + 1) + 0.076) \rceil \leq \lceil k \cdot (\log_3 d + 0.076) \rceil.$$

Now the proof of Proposition 2 has been finished. □

In the new version, the result $g_1(76|_{22}) = 87$ has been improved to $g_1(76|_{19}) = 75$ with new sub-algorithms $g_1(2,19,19) = 6$ and $g_1(4,4,4,4,76) = 9$.

**Acknowledgments.**
The author should be indebted to Professor L.K. Hua, Professor M. Aigner and Professor Y.P. Liu. Professor L.K. Hua was one of mathematical masters in Chinese modern mathematics, it was from his a popular lecture to the mid-school students I had known the counterfeit coins problem first time when I was a mid-school student. I was much benefited from Professor M. Aigner's nice book "combinatorial search ". Most of the researches above had been made in the period of my Ph.D. studies, Professor Y.P. Liu, the advisor for my Ph.D., gave me much supports and encouragements.

References

[1] A.P. Li, Some results on Counterfeit Coins Problem, arXiv, e-Print archive,0903.0905

# Appendix     Algorithms

## Algorithm  $g_1(5,4,4) = 4$

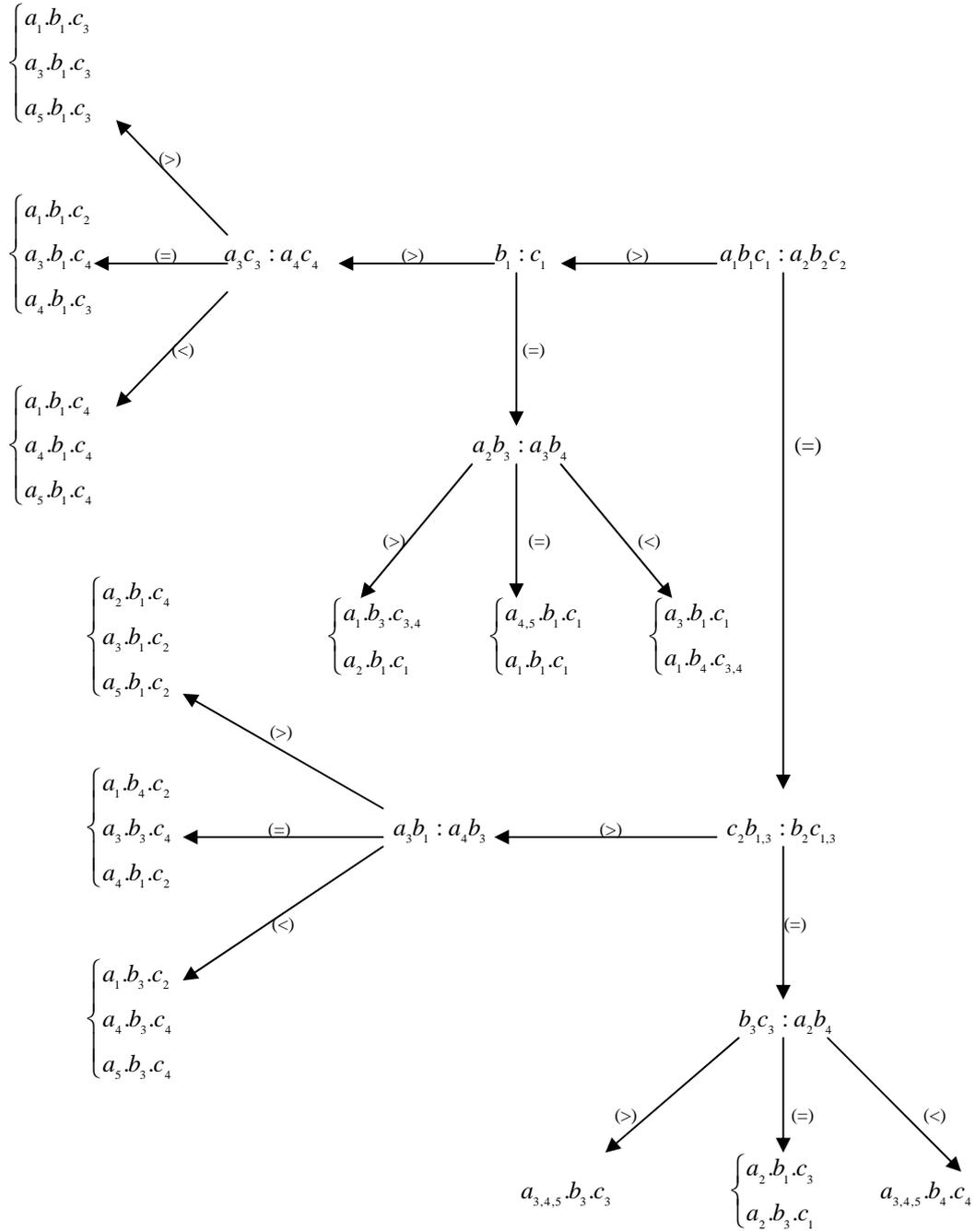

$$A = \{a_i\}_1^5, \quad B = \{b_i\}_1^4, \quad C = \{c_i\}_1^4, \quad Ct(A) = Ct(B) = Ct(C) = 1.$$

Fig. 1

Algorithm $(4|_5) \xrightarrow{5} (5)$

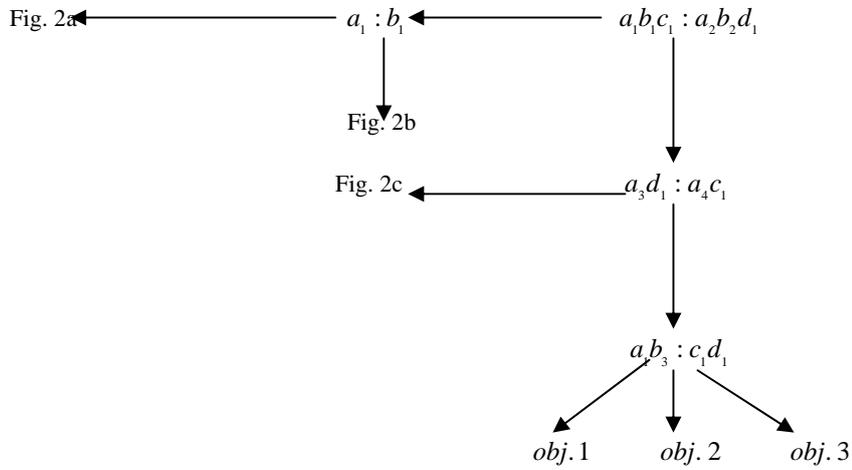

Fig. 2

$obj.1 = a_1.b_2.c_{2-4}.d_{2-4}.e_{1-4}, \quad obj.2 = a_2.b_1.c_{2-4}.d_{2-4}.e_{1-4}$

$obj.3 = a_1.b_2.c_1.d_1.e_{1-4} \cup a_2.b_1.c_1.d_1.e_{1-4} \cup a_3.b_2.c_1.d_{2-4}.e_{1-4} \cup a_4.b_1.c_{2-4}.d_1.e_{1-4}$

List 2

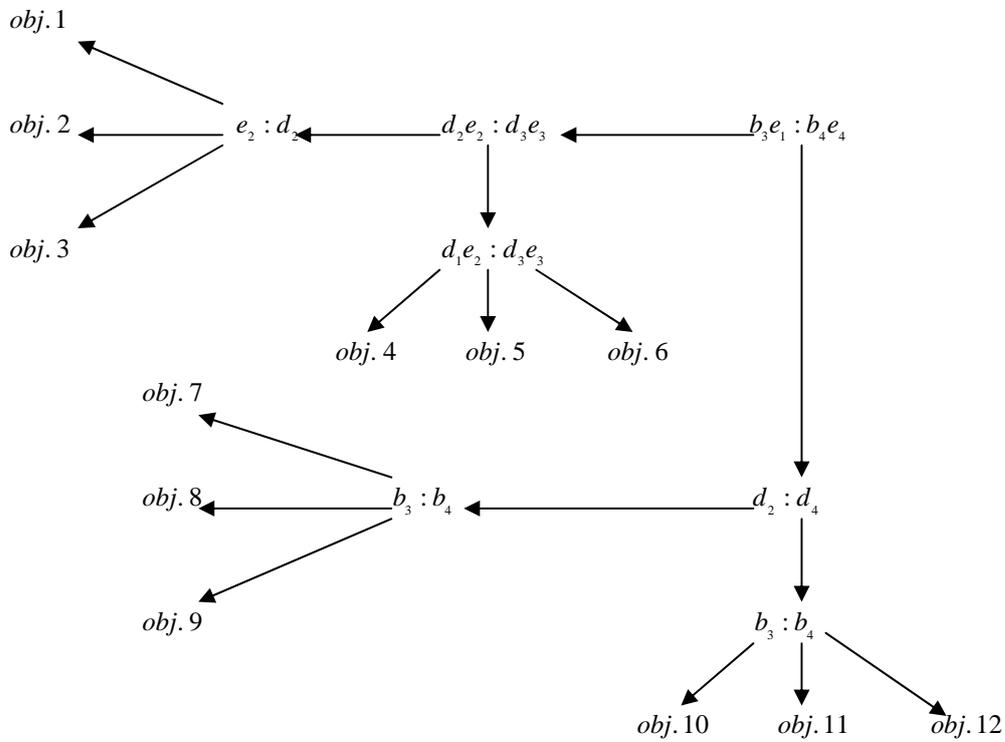

Fig. 2a

$$obj.1 = \begin{cases} a_1.b_3.c_{1-4}.d_4.e_2 \\ a_1.b_3.c_1.d_1.e_2 \end{cases}, \quad obj.2 = a_1.b_3.c_{1-4}.d_2.e_2, \quad obj.3 = \begin{cases} a_1.b_3.c_{1-4}.d_4.e_2 \\ a_1.b_3.c_1.d_1.e_2 \end{cases}$$

$$obj.4 = \begin{cases} a_1.b_3.c_{1-4}.d_4.e_2 \\ a_1.b_3.c_1.d_1.e_2 \end{cases}, \quad obj.5 = a_1.b_3.c_{1-4}.d_4.e_1, \quad obj.6 = \begin{cases} a_1.b_3.c_{1-4}.d_4.e_2 \\ a_1.b_3.c_1.d_1.e_2 \end{cases}$$

$$obj.7 = a_1.b_3.c_{1-4}.d_2.e_4, \quad obj.8 = a_1.b_2.c_1.d_2.e_{2,3}, \quad obj.9 = a_1.b_4.c_{1-4}.d_2.e_3,$$

$$obj.10 = \begin{cases} a_1.b_3.c_{1-4}.d_3.e_4 \\ a_1.b_3.c_1.d_1.e_4 \end{cases}, \quad obj.11 = a_1.b_2.c_1.d_3.e_{2,3}, \quad obj.12 = \begin{cases} a_1.b_4.c_{1-4}.d_3.e_3 \\ a_1.b_4.c_1.d_1.e_3 \end{cases}$$

List 2a

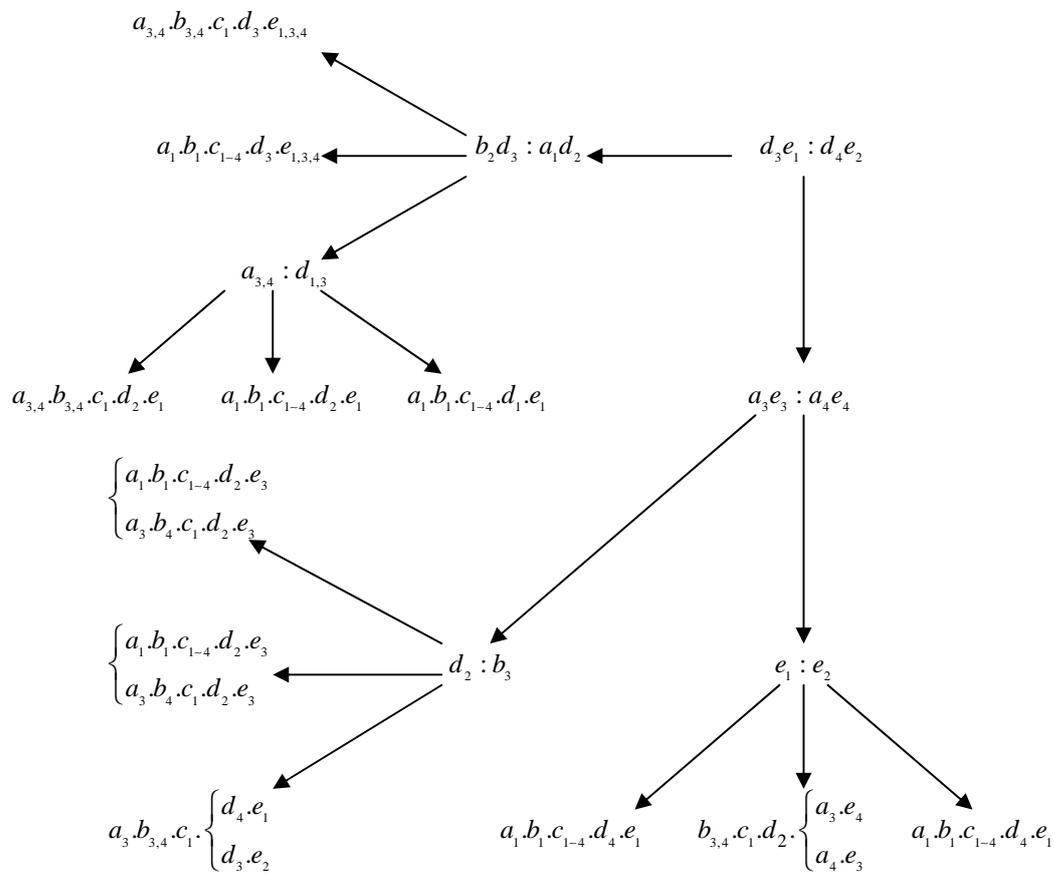

Fig. 2b

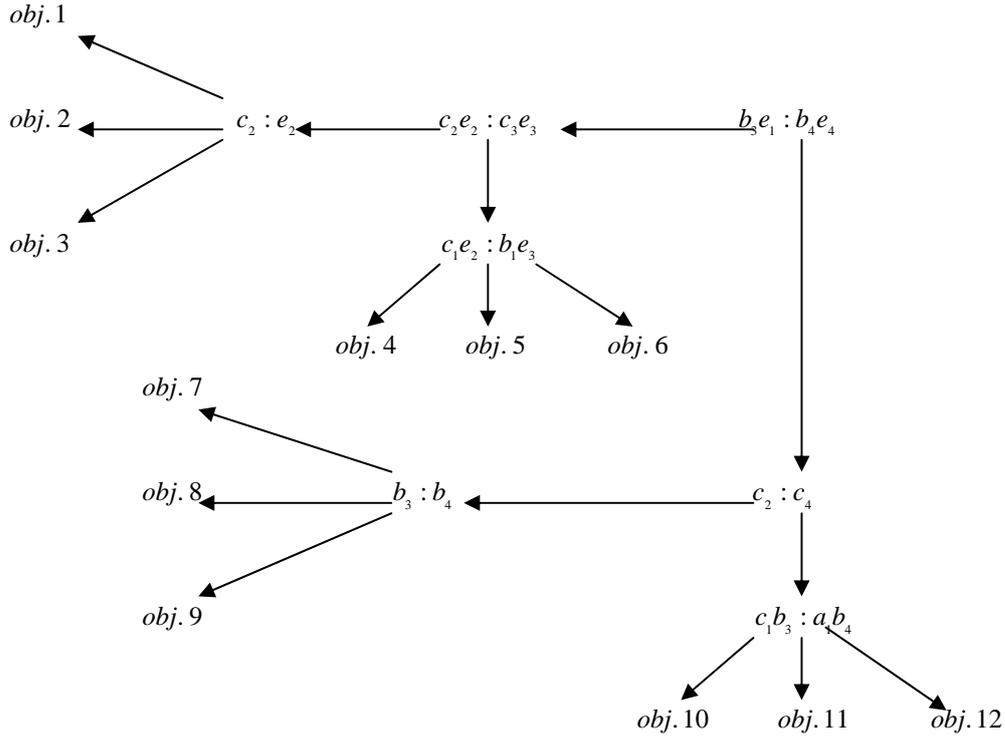

Fig. 2c

$$obj.1 = \begin{cases} a_1.b_3.c_2.d_1.e_4 \\ a_3.b_1.c_2.d_1.e_4 \\ a_3.b_3.c_2.d_{2-4}.e_1 \end{cases}, \quad obj.2 = \begin{cases} a_1.b_3.c_2.d_1.e_2 \\ a_3.b_3.c_2.d_{2-4}.e_2 \end{cases}, \quad obj.3 = \begin{cases} a_1.b_3.c_4.d_1.e_2 \\ a_3.b_3.c_1.d_1.e_2 \\ a_3.b_3.c_4.d_{2-4}.e_2 \end{cases}$$

$$obj.4 = \begin{cases} a_1.b_3.c_3.d_{2-4}.e_2 \\ a_1.b_3.c_3.d_1.e_2 \\ a_3.b_3.c_1.d_1.e_1 \end{cases}, \quad obj.5 = \begin{cases} a_1.b_3.c_4.d_1.e_1 \\ a_3.b_3.c_4.d_{2-4}.e_1 \end{cases}, \quad obj.6 = \begin{cases} a_3.b_3.c_2.d_{2-4}.e_3 \\ a_1.b_3.c_2.d_1.e_3 \\ a_3.b_1.c_4.d_1.e_1 \end{cases}$$

$$obj.7 = \begin{cases} a_3.b_3.c_2.d_{2-4}.e_4 \\ a_1.b_3.c_2.d_1.e_4 \end{cases}, \quad obj.8 = a_3.b_1.c_2.d_1.e_{2,3}, \quad obj.9 = \begin{cases} a_3.b_4.c_2.d_{2-4}.e_1 \\ a_1.b_4.c_2.d_1.e_1 \end{cases}$$

$$obj.10 = \begin{cases} a_3.b_3.c_3.d_{2-4}.e_4 \\ a_3.b_3.c_1.d_1.e_4 \end{cases}, \quad obj.11 = \begin{cases} a_1.b_3.c_3.d_1.e_4 \\ a_3.b_1.c_2.d_1.e_{2,3} \\ a_3.b_4.c_1.d_1.e_1 \end{cases}, \quad obj.12 = \begin{cases} a_3.b_4.c_3.d_{2-4}.e_1 \\ a_1.b_4.c_3.d_1.e_1 \end{cases}$$

List 2c

Algorithm $g_1(8,8,8,4) = 7$

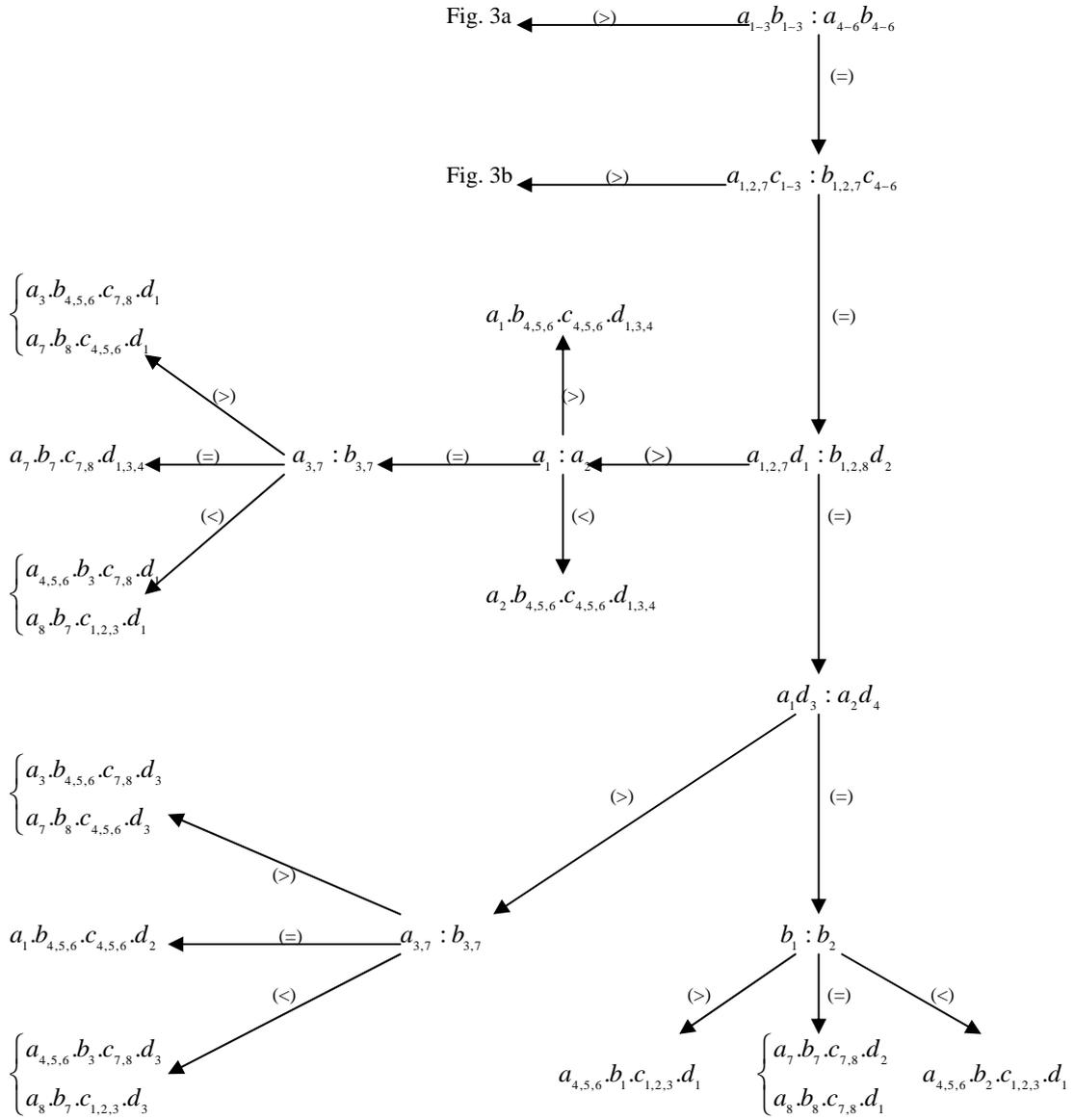

$A = \{a_i\}_1^8$, $B = \{b_i\}_1^8$, $C = \{c_i\}_1^8$, $D = \{d_i\}_1^4$, $Ct(A) = Ct(B) = Ct(C) = Ct(D) = 1$.

Fig. 3

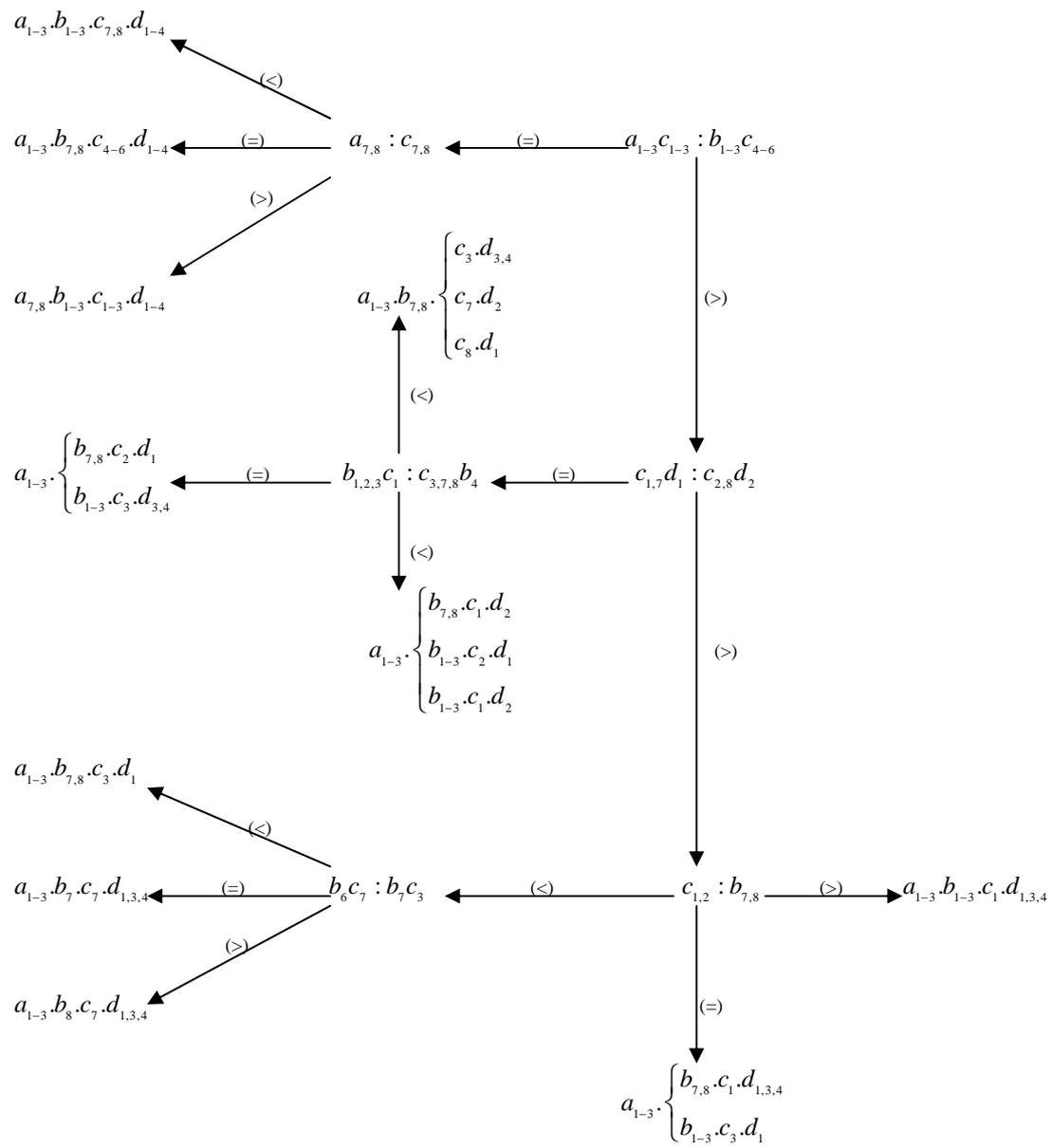

Fig. 3a

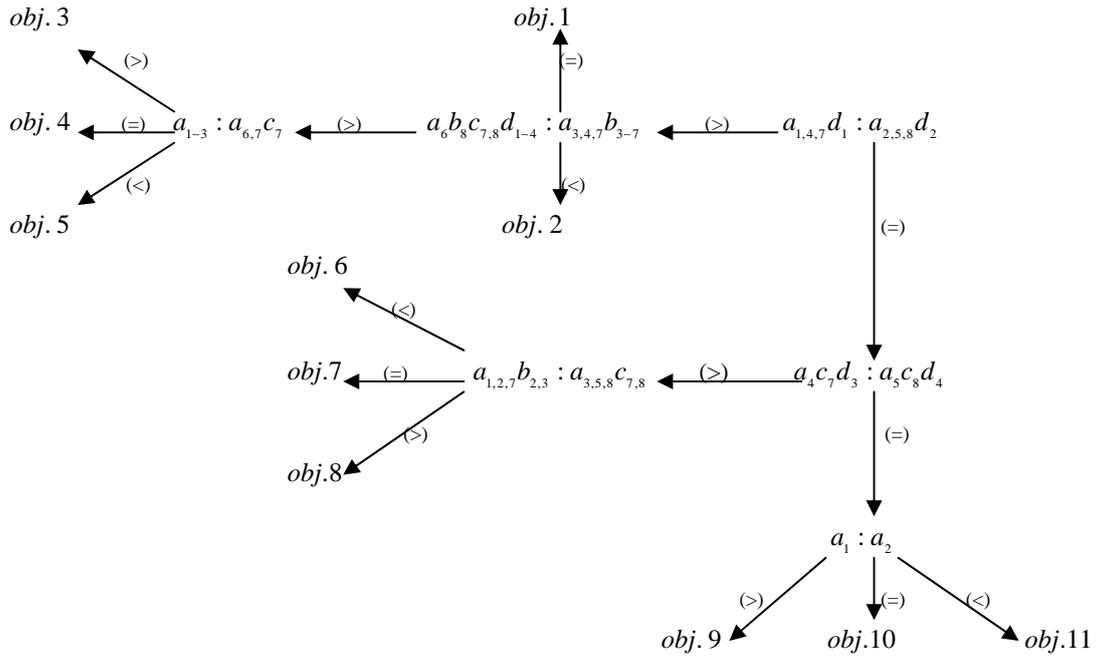

Fig. 3b

$obj.1 = a_1.b_{4-6}.c_{1-3}.d_{1,3,4}$, $\quad obj.2 = \begin{cases} a_3.b_{4-6}.c_{1-3}.d_1 \\ a_4.b_3.c_{1-3}.d_{1,3,4} \\ a_7.b_7.c_{1-3}.d_{1,3,4} \end{cases}$

$obj.3 = a_1.b_{4-6}.c_8.d_{1,3,4}$, $\quad obj.4 = a_1.b_{4-6}.c_7.d_{1,3,4}$, $\quad obj.5 = \begin{cases} a_7.b_8.c_{1,2,3,7,8}.d_1 \\ a_6.b_3.c_{1-3}.d_1 \end{cases}$

$obj.6 = c_{1-3}.\begin{cases} a_4.b_2.d_2 \\ a_6.b_3.d_3 \\ a_7.b_8.d_3 \end{cases}$, $\quad obj.7 = b_{4-6}.\begin{cases} a_2.c_7.d_1 \\ a_1.c_7.d_2 \end{cases}$, $\quad obj.8 = a_3.b_{4-6}.c_{1-3}.d_3$,

$obj.9 = a_1.b_{4-6}.c_{1-3}.d_2$, $\quad obj.10 = \begin{cases} c_{1-3}.a_7.b_7.d_2 \\ c_{1-3}.a_8.b_8.d_1 \\ a_7.b_8.c_8.d_3 \\ a_7.b_8.c_7.d_4 \end{cases}$, $\quad obj.11 = a_2.b_{4-6}.c_{1-3}.d_1$

List 3b

Algorithm $g_1(2,4,10) = 4$

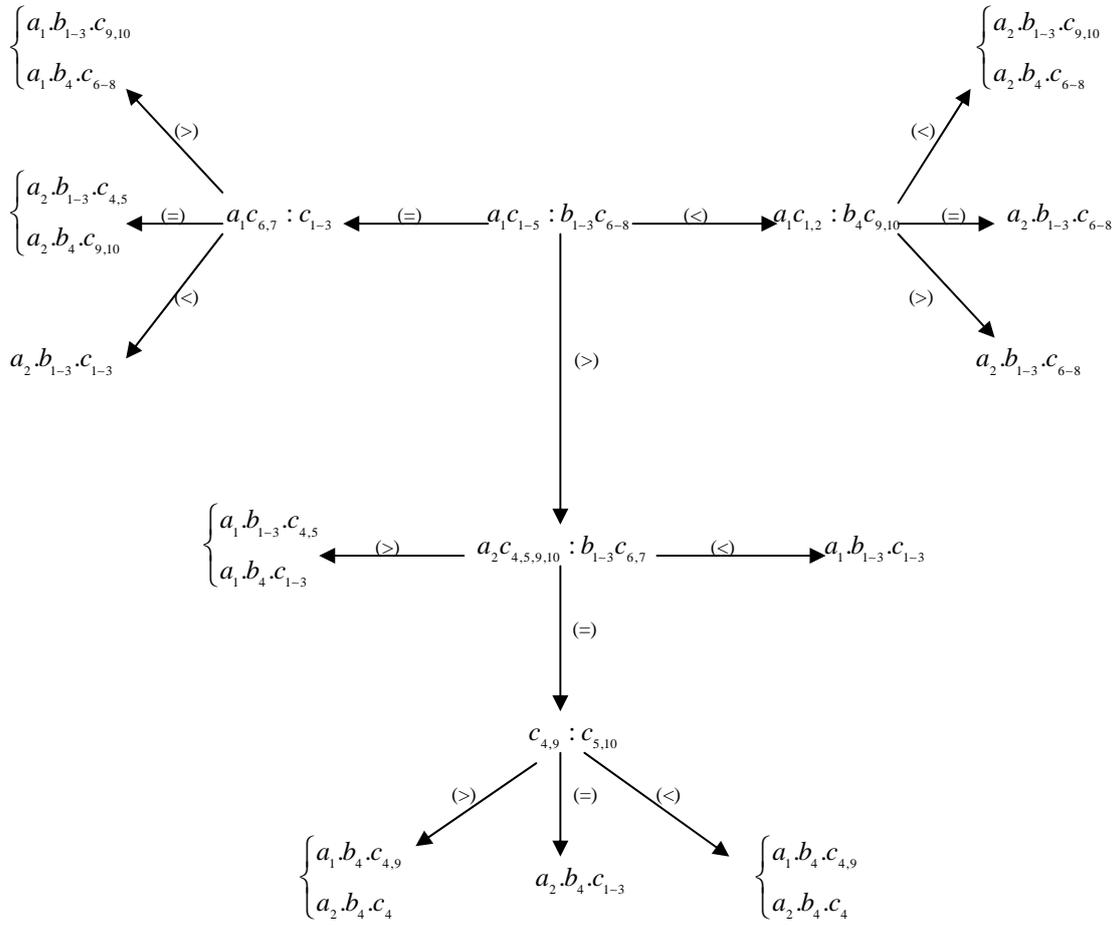

$$A = \{a_i\}_1^2,\ B = \{b_i\}_1^4,\ C = \{c_i\}_1^{10},\ Ct(A) = Ct(B) = Ct(C) = 1$$

Fig. 4

Algorithm $(16,16) \xrightarrow{3} (10)$

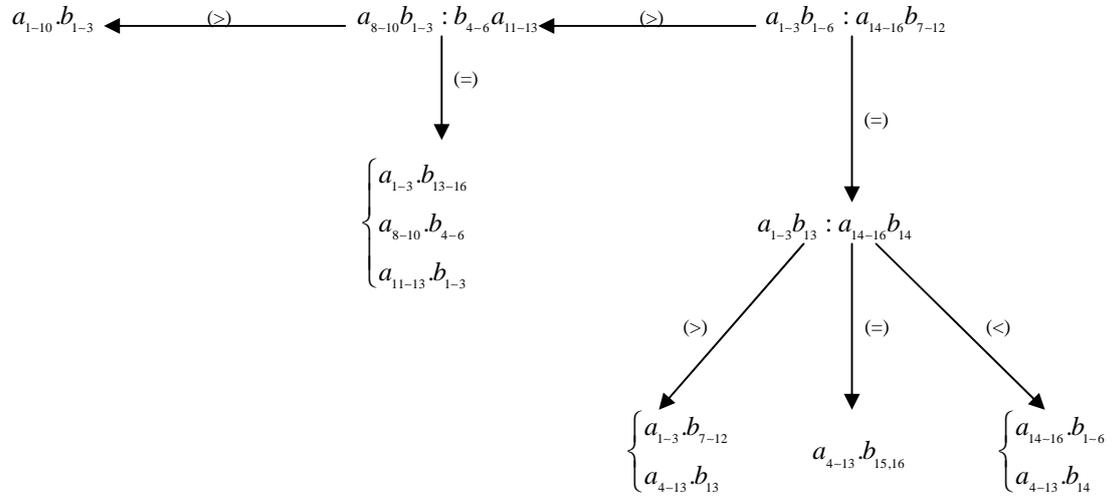

Fig. 5

Algorithm $(19,19) \xrightarrow{3} (14) \vee (2 \times 7)$

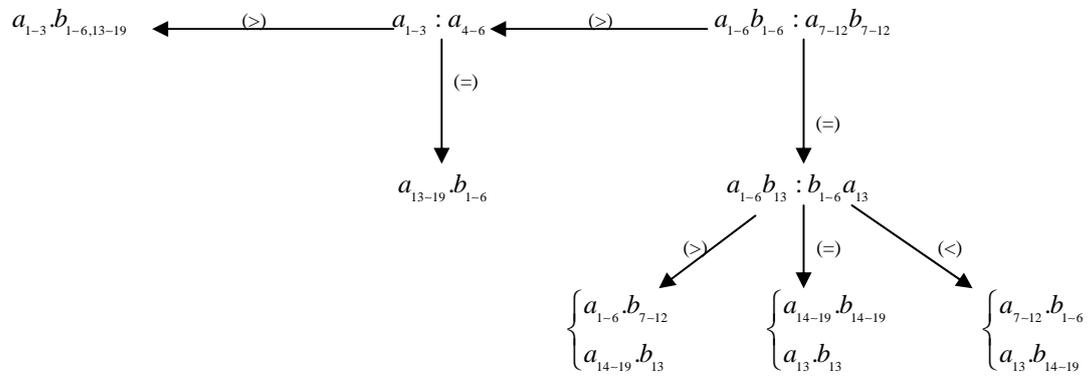

Fig. 6

Algorithm $g_1(11, 22) = 5$

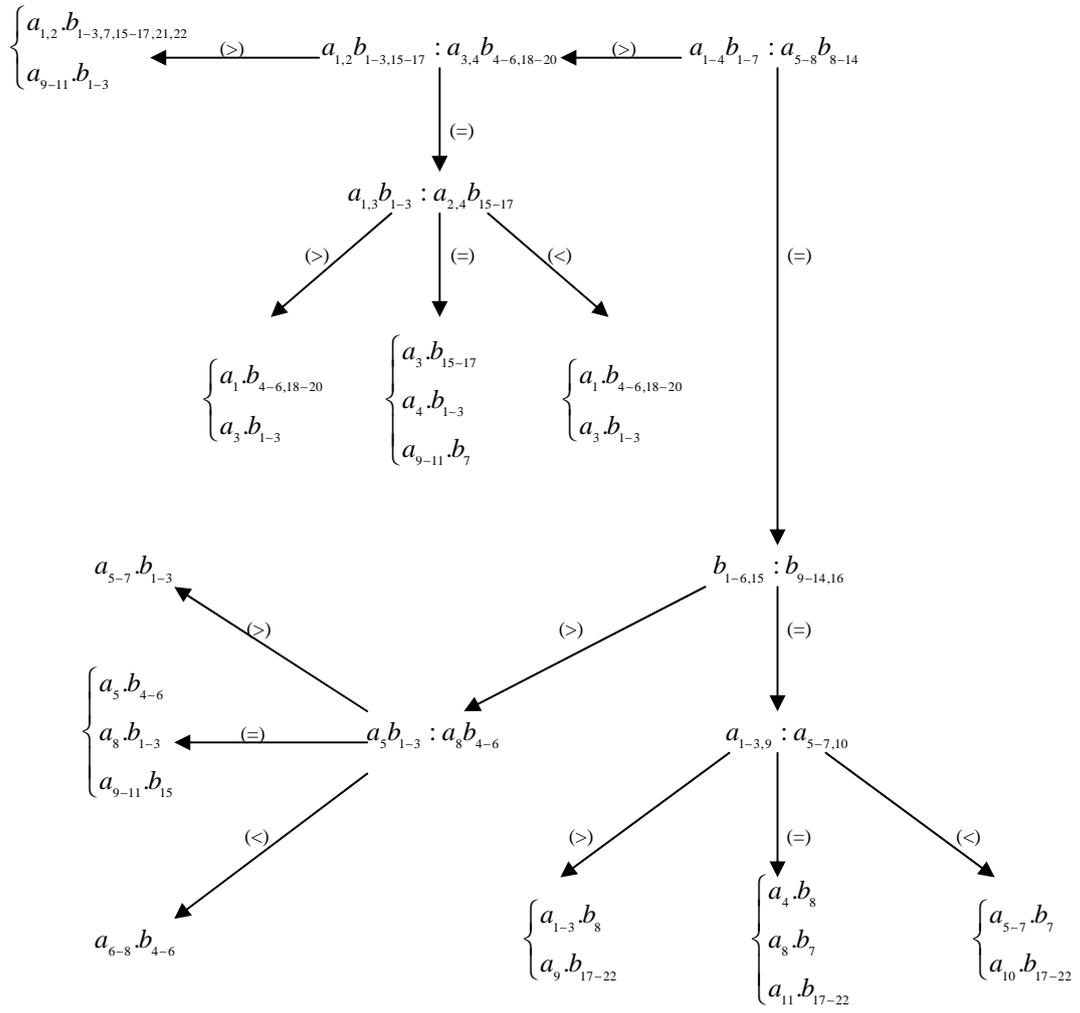

Fig. 7

## Sub-algorithm $g_1(4,4,4,4,76) = 9$

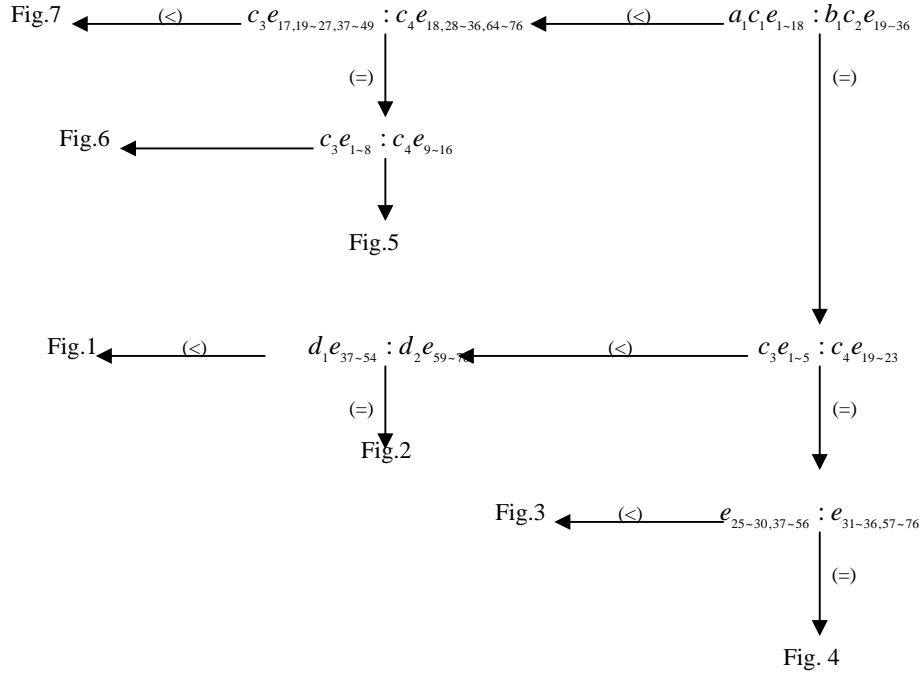

Fig. Main

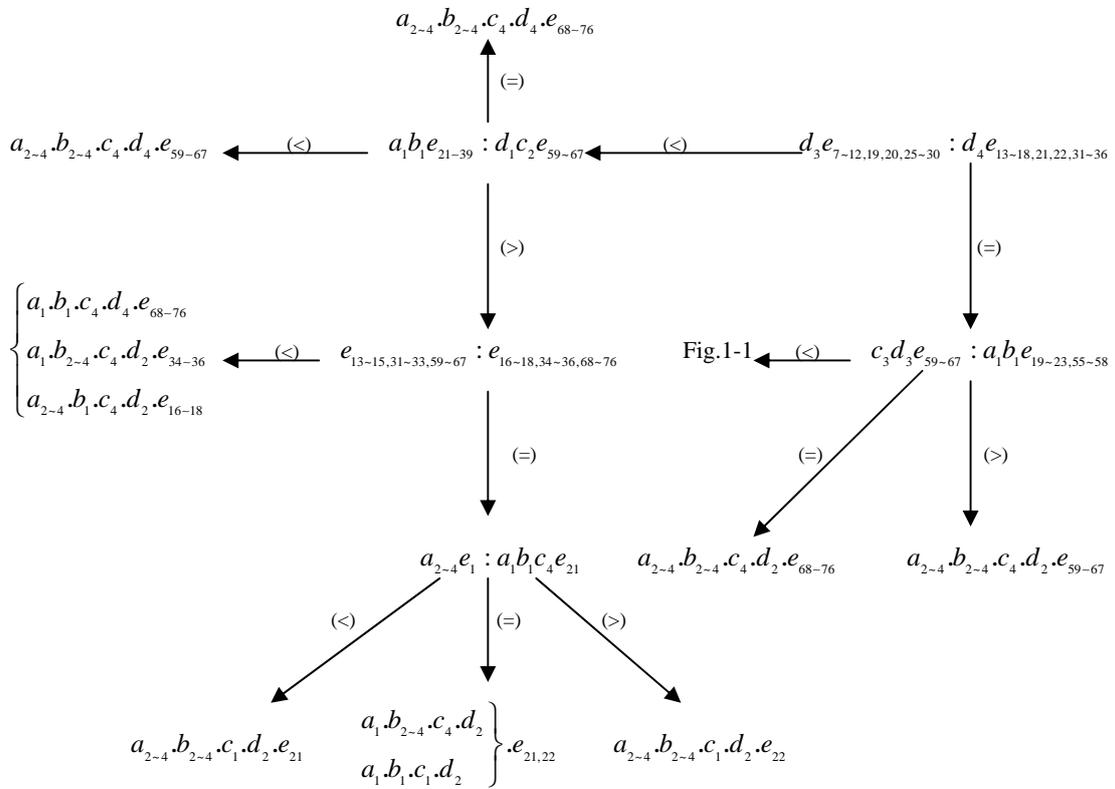

Fig. 1

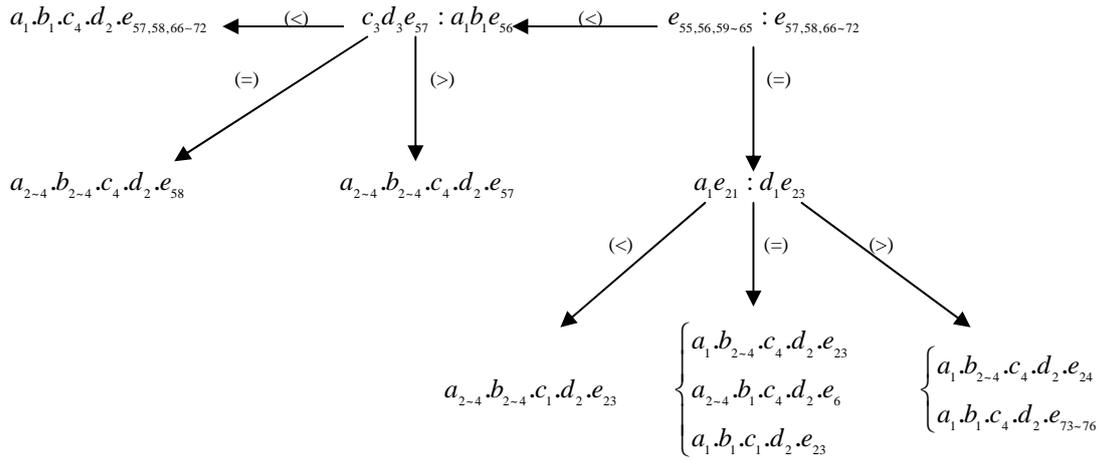

Fig. 1-1

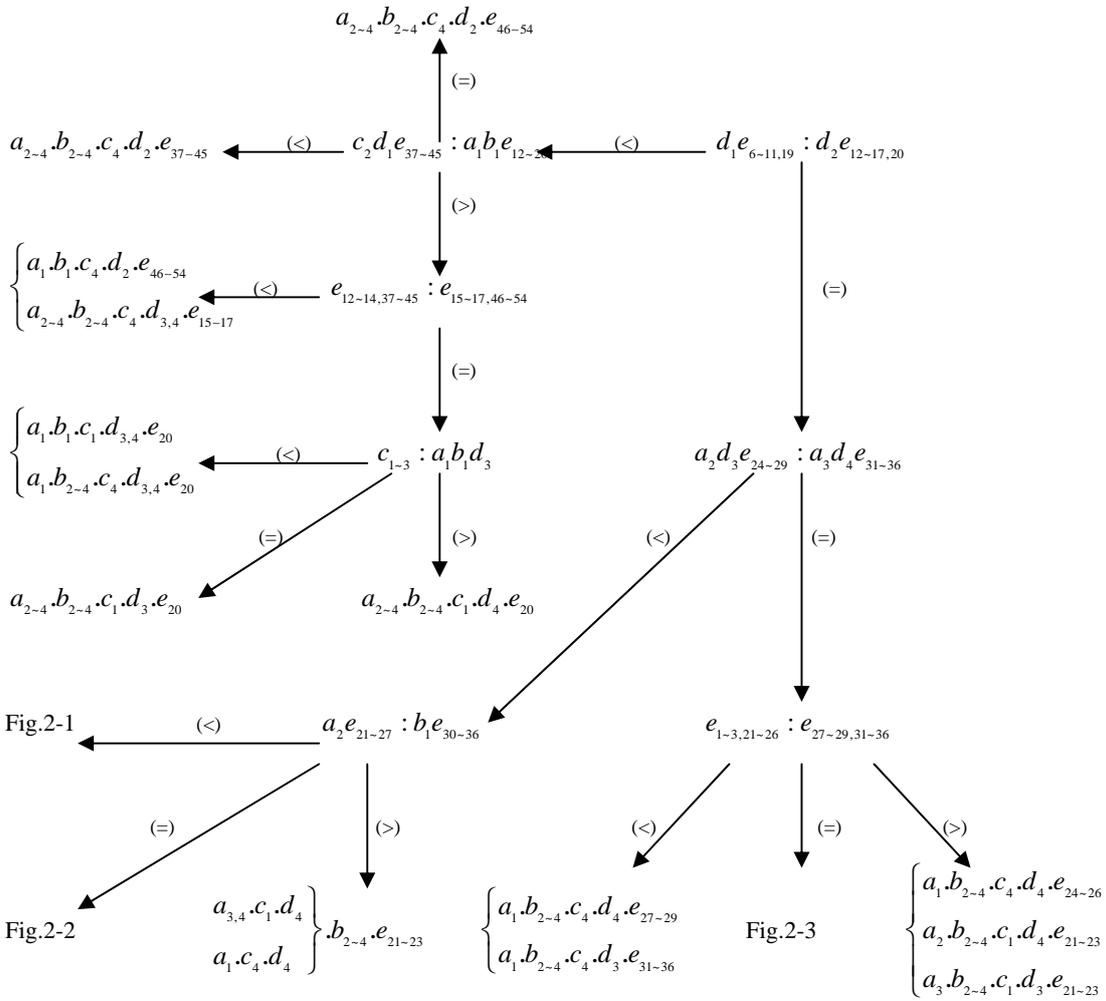

Fig.2

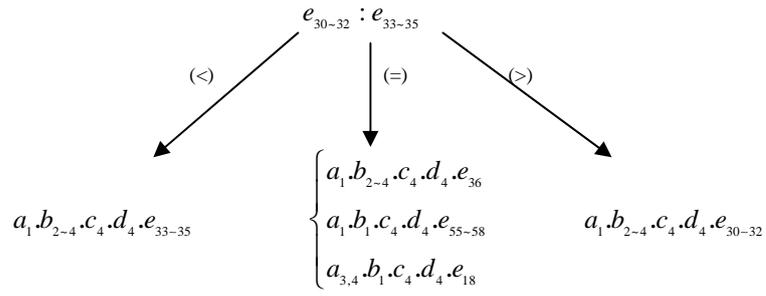

Fig.2-1

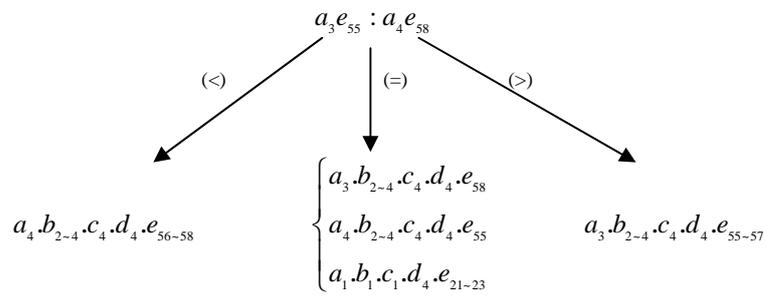

Fig.2-2

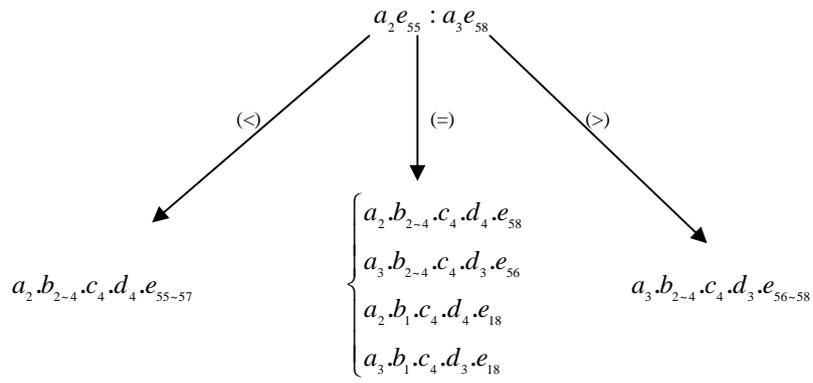

Fig.2-3

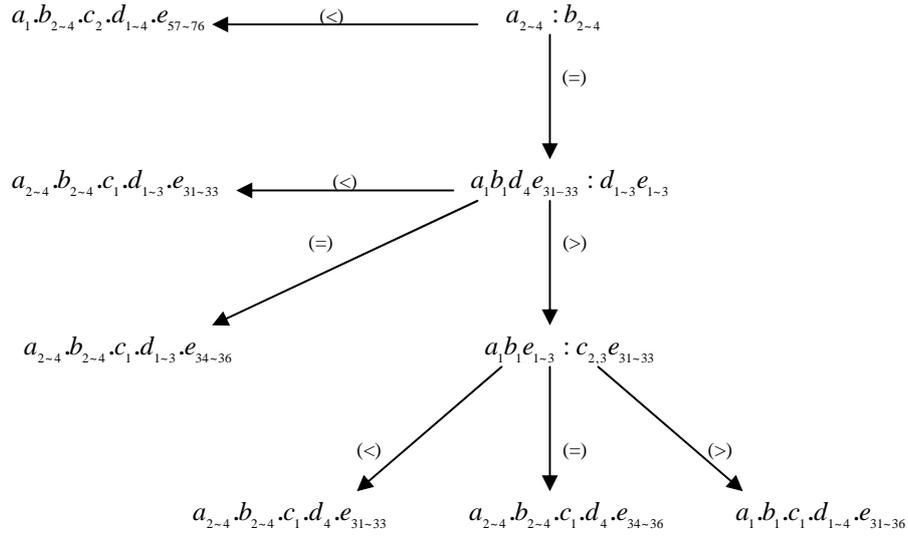

Fig. 3

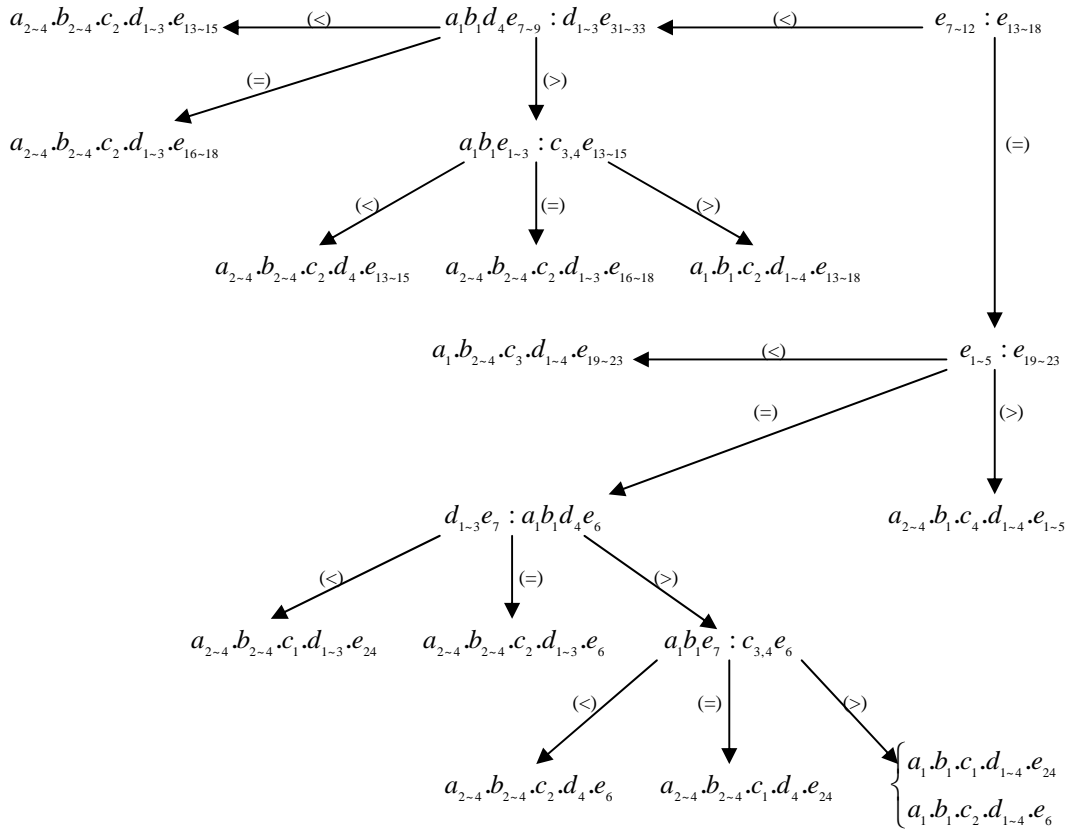

Fig.4

$a_{2\sim4}.b_{2\sim4}.c_2.d_{1\sim3}.e_{55\sim63}$

$\uparrow (>)$

$a_{1\sim4}.b_1.c_2.d_{1\sim4}.e_{50\sim54} \xleftarrow{(<)} d_4 e_{1\sim5} : b_1 e_{50\sim54} \xleftarrow{(<)} d_{1\sim3} e_{1\sim5} : c_1 b_1 d_4 e_{50\sim54}$

$a_{2\sim4}.b_{2\sim4}.c_2.d_4.e_{50\sim52}$

$\uparrow (>) \qquad (=) \downarrow \qquad (>) \downarrow \qquad (=) \downarrow$

$\begin{cases} a_1.b_1.c_2.d_4.e_{55\sim63} \\ a_{2\sim4}.b_{2\sim4}.c_2.d_4.e_{53,54} \end{cases} \xleftarrow{(=)} a_1 e_{50\sim52} : b_1 e_{1\sim3} \qquad a_{2\sim4}.b_{2\sim4}.c_2.d_4.e_{55\sim63}$

$\swarrow (<)$

$a_{2\sim4}.b_1.c_2.d_4.e_{55\sim63} \qquad\qquad a_{2\sim4}.b_1.c_2.d_{1\sim3}.e_{55\sim63} \xleftarrow{(<)} a_1 e_{50\sim52} : b_1 e_{1\sim3}$

$\qquad\qquad\qquad (=) \swarrow \qquad (>) \downarrow$

$\begin{cases} a_1.b_1.c_2.d_{1\sim3}.e_{55\sim63} \\ a_{2\sim4}.b_{2\sim4}.c_2.d_{1\sim3}.e_{53,54} \end{cases} \qquad a_{2\sim4}.b_{2\sim4}.c_2.d_{1\sim3}.e_{50\sim52}$

Fig. 5

$a_{2\sim4}.b_{2\sim4}.c_4.d_4.e_{19\sim27} \qquad\qquad a_{2\sim4}.b_{2\sim4}.c_4.d_{1\sim3}.e_{19\sim27}$

$\uparrow (>) \qquad\qquad\qquad\qquad \uparrow (>)$

$\begin{cases} a_{2\sim4}.b_1.c_4.d_{1\sim3}.e_{37\sim43} \\ a_{2\sim4}.b_1.c_4.d_4.e_{44\sim49} \end{cases} \xleftarrow{(<)} b_{2\sim4} e_{1\sim6} : d_{1\sim3} e_{44\sim49} \xleftarrow{(<)} d_{1\sim3} e_{1\sim7} : c_1 b_1 d_4 e_{37\sim43}$

$(=) \downarrow$

$a_{2\sim4}.b_1.c_4.d_4.e_{19\sim27} \xleftarrow{(<)} a_1 e_{9\sim16,37} : b_2 e_{19\sim27} \qquad\qquad (=) \downarrow$

$(=) \swarrow \qquad (>) \downarrow$

$\begin{cases} a_1.b_1.c_4.d_4.e_{19\sim27} \\ a_{2\sim4}.b_1.c_4.d_4.e_{38\sim43} \end{cases} \qquad a_{2\sim4}.b_1.d_4. \begin{cases} c_2.e_{9\sim16} \\ c_4.e_{37} \end{cases}$

$a_{2\sim4}.b_1.c_4.d_{1\sim3}.e_{19\sim27} \xleftarrow{(<)} a_1 e_{41\sim49} : b_2 e_{19\sim27}$

$(=) \swarrow \qquad (>) \downarrow$

$\begin{cases} a_1.b_1.c_4.d_{1\sim3}.e_{19\sim27} \\ a_{2\sim4}.b_1.c_4.d_{1\sim3}.e_{44\sim49} \end{cases} \qquad a_{2\sim4}.b_1.c_2.d_{1\sim3}.e_{9\sim16}$

Fig.6

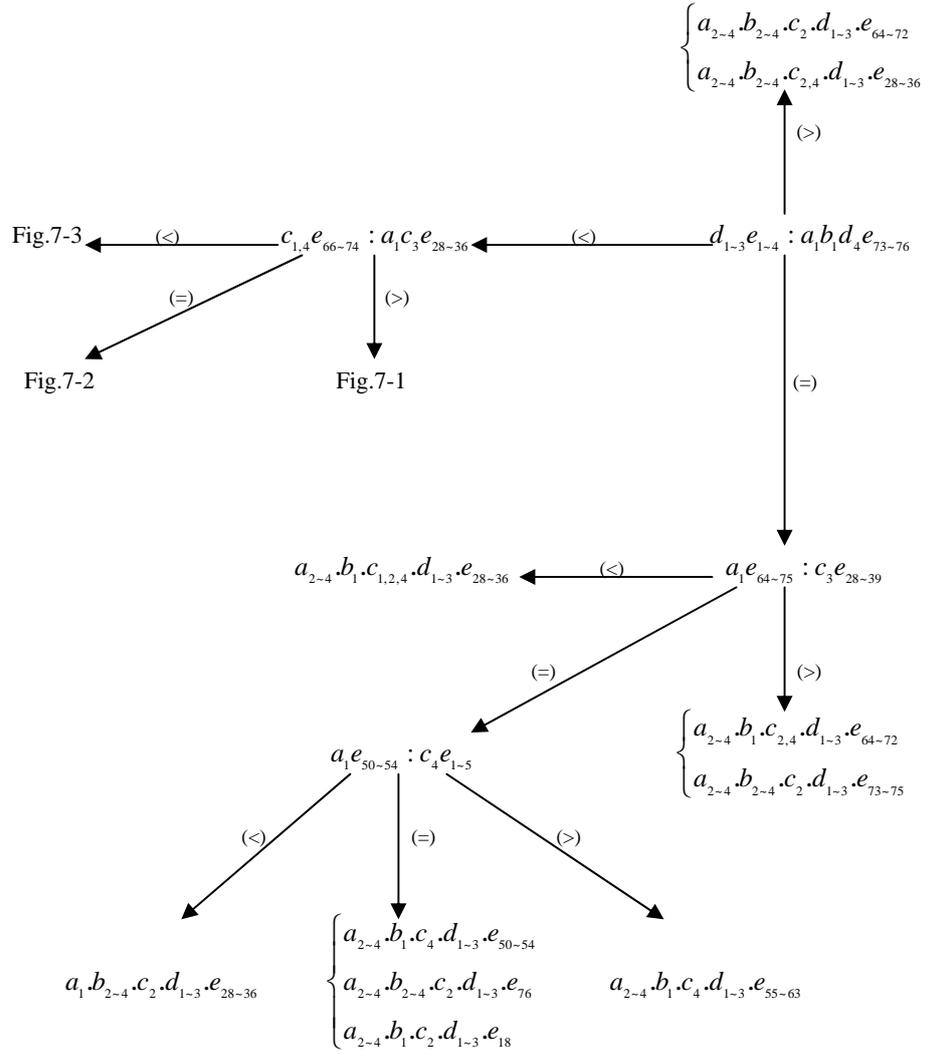

Fig.7

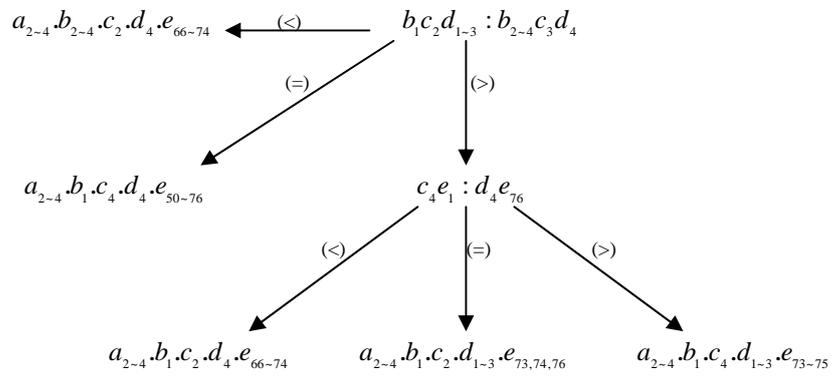

Fig.7-1

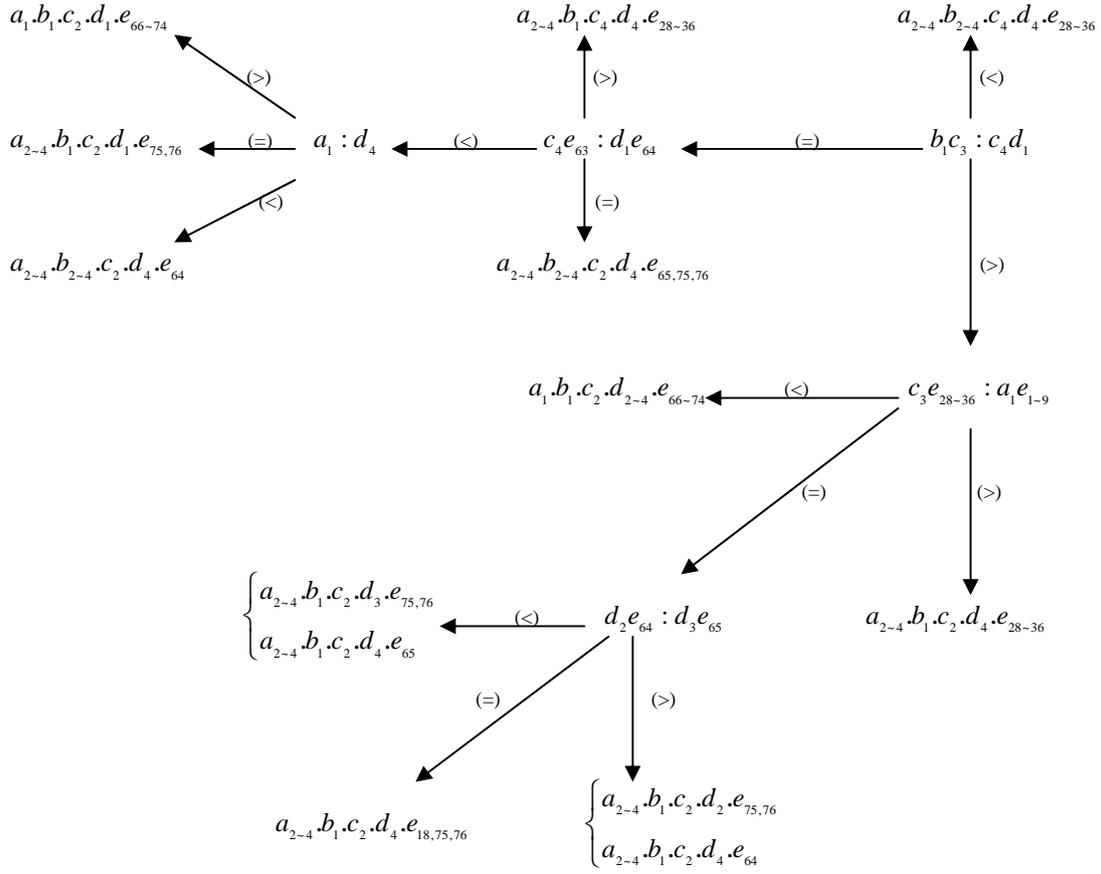

Fig.7-2

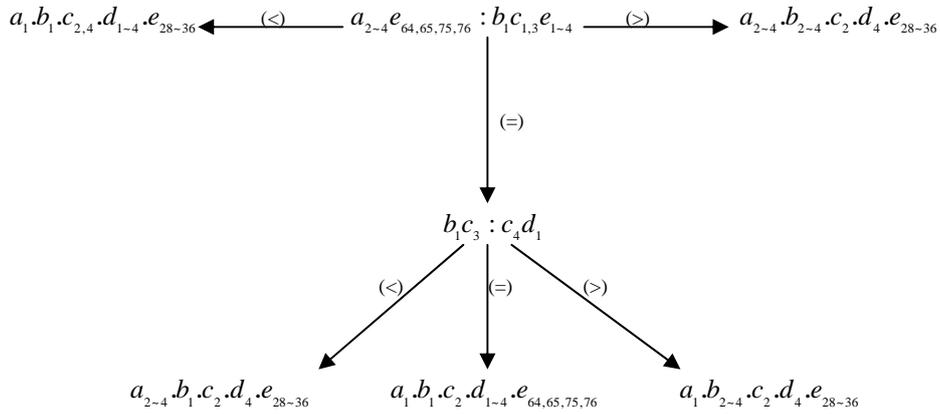

Fig. 7-3

$A = \{a_i\}_1^4,\ B = \{b_i\}_1^4,\ C = \{c_i\}_1^4,\ D = \{d_i\}_1^4,\ E = \{e_i\}_1^{76},$
$Ct(A) = Ct(B) = Ct(C) = Ct(D) = Ct(E) = 1.$

Sub-algorithm $g_1(2,19,19) = 6$

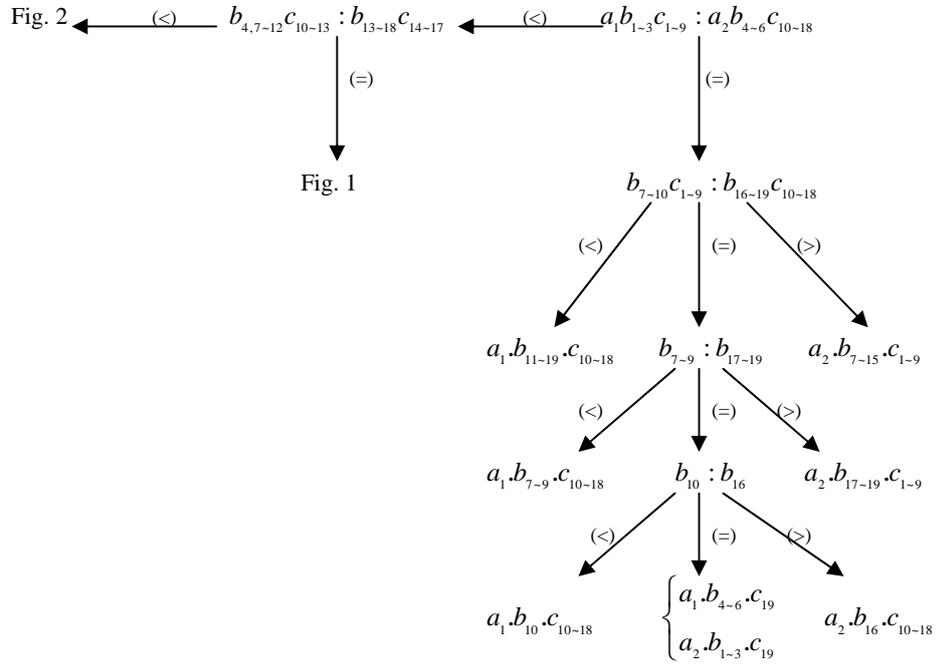

Fig. Main

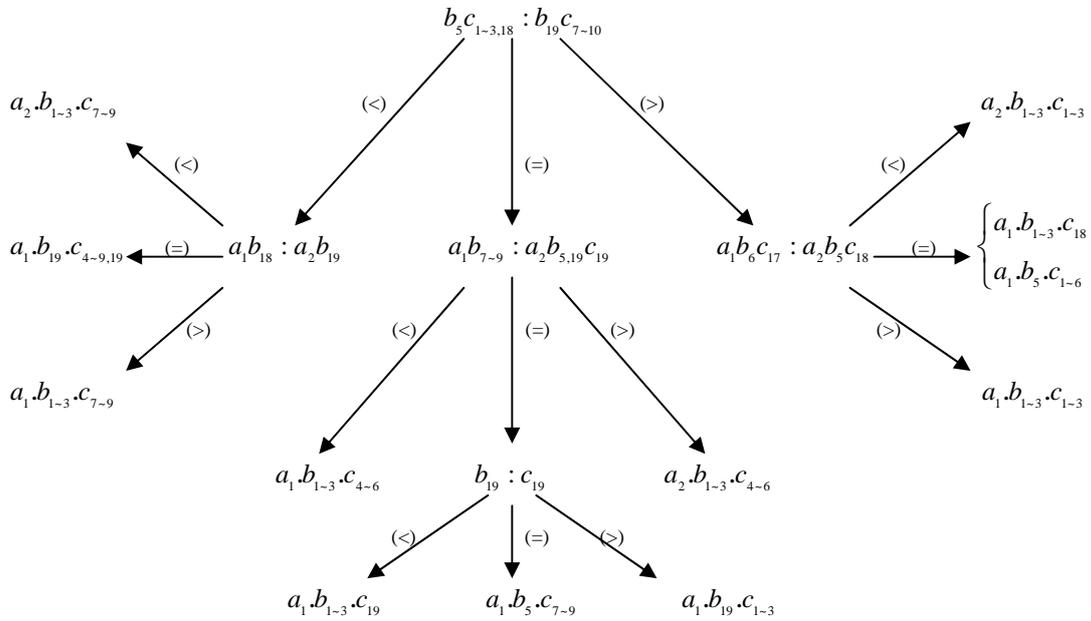

Fig.1

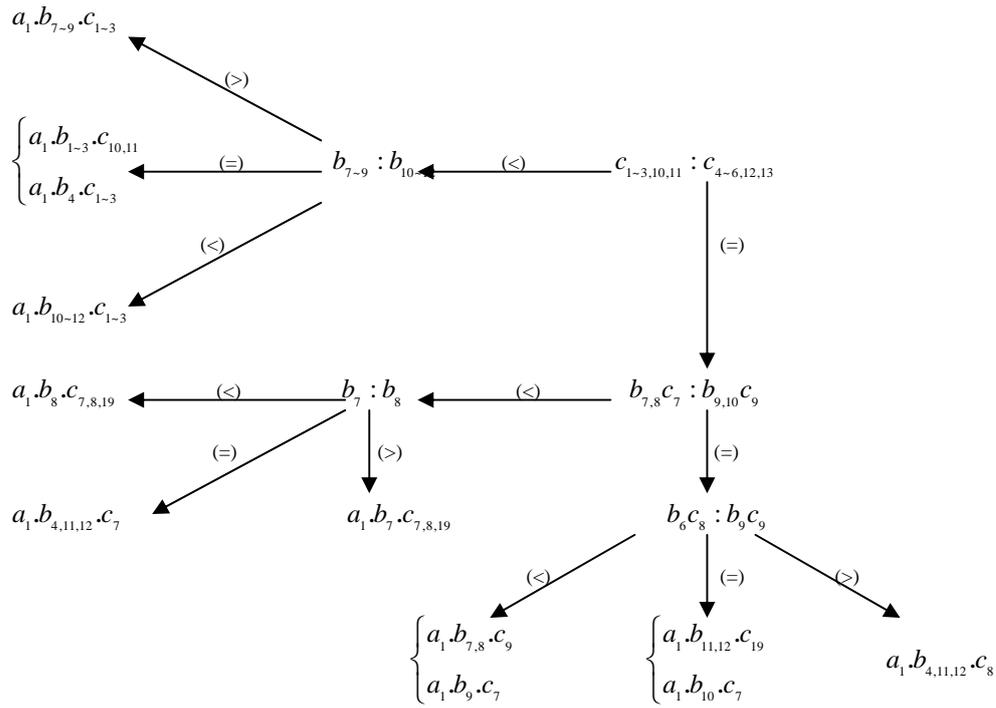

Fig. 2

$A = \{a_i\}_1^2$, $B = \{b_i\}_1^{19}$, $C = \{c_i\}_1^{19}$,
$Ct(A) = Ct(B) = Ct(C) = 1$.

To save the space, some of algorithm schemas have been omitted.